\documentclass[a4paper,12pt]{article}
\usepackage{amsmath}
\usepackage{amsmath,amssymb,latexsym,color}
\usepackage[mathscr]{eucal}

\newtheorem{thm}{Theorem}{}

\newtheorem{lemma}[thm]{Lemma}
\newtheorem{cor}[thm]{Corollary}

\newtheorem{example}{Example}[section]
\newtheorem{defin}{Definition}[section]

\newcommand{\proof}{{\it Proof.\quad}}
\newcommand{\qed}{\hfill\Box\medskip}

\renewcommand{\abovewithdelims}[2]{%
\genfrac{[}{]}{0pt}{}{#1}{#2}}

\usepackage{CJK}

\begin{document}
\begin{CJK*}{GBK}{song}

\renewcommand{\baselinestretch}{1.3}
\title{\bf The eigenvalues of $q$-Kneser graphs}

\author{
Benjian Lv\quad Kaishun Wang\footnote{Corresponding author. E-mail address: wangks@bnu.edu.cn}\\
{\footnotesize   \em  Sch. Math. Sci. {\rm \&} Lab. Math. Com. Sys.,
Beijing Normal University, Beijing, 100875,  China} }
 \date{}
 \maketitle

\begin{abstract}
In this note, we prove some combinatorial identities and obtain a
simple form of the eigenvalues of $q$-Kneser graphs.



\end{abstract}

Let $q$ be a prime power. For any integer $n$ and positive integer
$i$, the Gaussian coefficients is defined by
$$
\abovewithdelims{n}{i}_q=\prod^{i-1}_{j=0}\frac{q^{n-j}-1}{q^{i-j}-1}.
$$
By convention $\abovewithdelims{n}{0}_q=1$ for every integer $n$.
From now on, we will omit the subscript $q$. Note that
\begin{eqnarray}\label{c1}
\abovewithdelims{n}{i}=\abovewithdelims{n-1}{i-1}+q^i\abovewithdelims{n-1}{i}.
\end{eqnarray}

Let $\mathbb{F}_q^v$ be a $v$-dimensional vector space over a finite
field $\mathbb{F}_q$. The \emph{$q$-Kneser graph} $qK(v,k)$
 has as vertex set  the collection of $k$-dimensional subspaces of $\mathbb{F}_q^v$. Two
vertices  are adjacent if they intersect trivially.  If $k\leq
v<2k$,  then $qK(n,k)$ is null graph, so we only consider the case
$v\geq 2k$.

Delsarte \cite{Delsarte} calculated the eigenvalues of Grassmann
schemes. In particular, the eigenvalues of $qK(v,k)$ were given.
\begin{thm}\label{get}{\rm(\cite[Theorem 10]{Delsarte})}
All the distinct eigenvalues of $qK(v,k)$ are
\begin{eqnarray}\label{c2}
\lambda_j=(-1)^jq^{(k-j)j+{j\choose
2}}\sum_{s=0}^{k-j}(-1)^sq^{s\choose2}\abovewithdelims{k-j}{s}\abovewithdelims{v-2j-s}{v-k-j},
\end{eqnarray}
where $j=0,1,\ldots,k$.
\end{thm}

The eigenvalues of Kneser graphs are deduced to a simple form in
\cite[Theorem 9.4.3]{Godsil}. In this note, we shall generalize this
result to vector spaces, and obtain a simple form of eigenvalues of
$qK(v,k)$. We start with some useful combinatorial identities.

\begin{lemma}\label{a1}
For any integer $n$ and nonnegative integer $i$, we have
$$
\abovewithdelims{n}{i}=(-1)^{i}q^{ni-{i\choose2}}\abovewithdelims{-n+i-1}{i}.
$$
\end{lemma}

\proof  If $i=0$, the identity is obvious. If $i>0$, then
\begin{eqnarray*}
\abovewithdelims{n}{i}&=&(-1)^{i}\prod^{i-1}_{j=0}\frac{1-q^{n-j}}{q^{i-j}-1}\\
&=&(-1)^{i}q^{\frac{(2n-i+1)i}{2}}\prod^{i-1}_{j=0}\frac{q^{-n+j}-1}{q^{i-j}-1}\\
&=&(-1)^{i}q^{ni-{i\choose2}}\prod^{i-1}_{j=0}\frac{q^{-n+i-1-j}-1}{q^{i-j}-1}\\
&=&(-1)^{i}q^{ni-{i\choose2}}\abovewithdelims{-n+i-1}{i},
\end{eqnarray*}
as desired. $\qed$

The following identity is a generalization of
\cite[Theorem~2.14]{Wan}.
\begin{lemma}\label{a2}
For any integer $n$ and nonnegative integer $a$, we have
\begin{eqnarray*}\label{eq:a2}
\sum^{a}_{s=0}(-1)^{s}q^{s\choose2}\abovewithdelims{n}{s}=q^{na}\abovewithdelims{a-n}{a}.
\end{eqnarray*}
\end{lemma}

\proof
 We prove the result by induction on $a$. If $a=0$, then
the result is trivial. Suppose $a\geq 1$. By induction and Lemma
\ref{a1},
\begin{eqnarray*}
\sum^{a}_{s=0}(-1)^{s}q^{s\choose2}\abovewithdelims{n}{s}
&=&\sum^{a-1}_{s=0}(-1)^{s}q^{s\choose2}\abovewithdelims{n}{s}+(-1)^{a}q^{a\choose2}\abovewithdelims{n}{a}\\
&=&q^{n(a-1)}\abovewithdelims{a-1-n}{a-1}+q^{na}\abovewithdelims{a-1-n}{a}\\
&=&q^{na}\abovewithdelims{a-n}{a}.
\end{eqnarray*}
Hence, the desired result follows. $\qed$

\begin{lemma}\label{b1}
Let $m$, $a$, $t$ be nonnegative integers with $t\leq a\leq m$. Then
\begin{equation}\label{eq:b1}
\sum_{s=0}^a(-1)^sq^{{s}\choose{2}}\abovewithdelims{m}{s}\abovewithdelims{a-s}{t}=q^{m(a-t)}\abovewithdelims{a-m}{a-t}.
\end{equation}
\end{lemma}

\proof
 We prove the result by induction on $a$ and $t$. If $t=0$,
(\ref{eq:b1}) is immediate by Lemma \ref{a2}. If $a=t$,
(\ref{eq:b1}) is straightforward. Now suppose $1\leq t<a$. By
(\ref{c1}) and induction,
\begin{eqnarray*}
& &\sum_{s=0}^a(-1)^sq^{{s}\choose{2}}\abovewithdelims{m}{s}\abovewithdelims{a-s}{t}\\
&=&\sum_{s=0}^{a-1}(-1)^sq^{{s}\choose{2}}\abovewithdelims{m}{s}\abovewithdelims{a-s}{t}\\
& =&\sum_{s=0}^{a-1}(-1)^sq^{{s}\choose{2}}\abovewithdelims{m}{s}\abovewithdelims{a-1-s}{t-1}+q^t\sum_{s=0}^{a-1}(-1)^sq^{{s}\choose{2}}\abovewithdelims{m}{s}\abovewithdelims{a-1-s}{t}\\
& = & q^{m(a-t)}\abovewithdelims{a-1-m}{a-t}+q^tq^{m(a-1-t)}\abovewithdelims{a-1-m}{a-1-t}\\
& = & q^{m(a-t)}\abovewithdelims{a-m}{a-t}.
\end{eqnarray*}
Hence, the desired result follows. $\qed$

\begin{thm}\label{b2}
Let $m$, $a$, $t$ be nonnegative integers with $a\geq m$ and $a\geq
t$. Then
\begin{equation}\label{eq:b2}
\sum_{s=0}^m(-1)^sq^{{s}\choose{2}}\abovewithdelims{m}{s}\abovewithdelims{a-s}{t}=q^{m(a-t)}\abovewithdelims{a-m}{a-t}.
\end{equation}
\end{thm}

\proof
 We prove the result by induction on $a$ and $t$. If $t=0$,
(\ref{eq:b2}) is immediate from Lemma \ref{a2}. If $a=m$,
(\ref{eq:b2}) holds by Lemma \ref{b1}. If $a=t$, (\ref{eq:b2}) is
straightforward. Suppose $a\geq m+1$ and $1\leq t<a$. By (\ref{c1})
and induction,
\begin{eqnarray*}
 &
&\sum_{s=0}^m(-1)^sq^{{s}\choose{2}}\abovewithdelims{m}{s}\abovewithdelims{a-s}{t}\\
& =&\sum_{s=0}^m(-1)^sq^{{s}\choose{2}}\abovewithdelims{m}{s}\abovewithdelims{a-1-s}{t-1}+q^t\sum_{s=0}^m(-1)^sq^{{s}\choose{2}}\abovewithdelims{m}{s}\abovewithdelims{a-1-s}{t}\\
& = & q^{m(a-t)}\abovewithdelims{a-1-m}{a-t}+q^tq^{m(a-1-t)}\abovewithdelims{a-1-m}{a-t-1}\\
& = & q^{m(a-t)}\abovewithdelims{a-m}{a-t}.
\end{eqnarray*}
Therefore, (\ref{eq:b2}) holds. $\qed$

Substituting  $t=a-m$ in (\ref{eq:b2}), we obtain
\begin{cor}\label{b3}
For nonnegative integers $a\geq m$, we have
\begin{equation*}\label{eq:b3}
\sum_{s=0}^m(-1)^sq^{{s}\choose{2}}\abovewithdelims{m}{s}\abovewithdelims{a-s}{a-m}=q^{m^2}\abovewithdelims{a-m}{m}.
\end{equation*}
\end{cor}

Next, we shall deduce a simple form of the eigenvalues $\lambda_j$
in (\ref {c2}) of $qK(v,k)$.
\begin{thm}\label{main} All the distinct eigenvalues of $qK(v,k)$
are
$$
\lambda_j=(-1)^jq^{{k\choose2}+{k-j+1\choose
2}}\abovewithdelims{v-k-j}{v-2k},\ \ \ j=0,1,\ldots,k.
$$
Moreover, the multiplicity of $\lambda_j$ is $1$ if $j=0$, and
$\abovewithdelims{v}{j}-\abovewithdelims{v}{j-1}$ if $j\geq 1$.
\end{thm}

\proof  By (\ref{c2}) and Corollary~\ref{b3}, we have
\begin{eqnarray*}
\lambda_j&=&(-1)^jq^{(k-j)j+{j\choose
2}}\sum_{s=0}^{k-j}(-1)^sq^{s\choose2}\abovewithdelims{k-j}{s}\abovewithdelims{v-2j-s}{v-2j-(k-j)}\\
&=&(-1)^jq^{(k-j)j+{j\choose2}}q^{(k-j)^2}\abovewithdelims{v-2j-(k-j)}{k-j}\\
 &=&(-1)^jq^{{k\choose2}+{k-j+1\choose2}}\abovewithdelims{v-k-j}{v-2k}.
\end{eqnarray*}

By similar arguments in \cite{Delsarte,Richard}, the multiplicity of
each $\lambda_j$ may be computed.$\qed$

\section*{Acknowledgement}
This research is partially supported by NCET-08-0052, NSF of China
(10871027) and the Fundamental Research Founds for the Central
Universities of China.

\end{CJK*}

\end{document}